\documentclass[12pt]{amsart}

\usepackage{amssymb}
\usepackage{enumerate}
\usepackage{amsthm}
\usepackage{amsmath}
\usepackage{mathtools}
\usepackage{txfonts}
\usepackage{longtable}
\usepackage{booktabs}

\newcommand{\av}{\check{a}}
\newcommand{\hv}{\check{h}}

\newcommand{\bC}{\mathbb{C}}
\newcommand{\bN}{\mathbb{N}}
\newcommand{\bP}{\mathbb{P}}

\newcommand{\bCx}{\bC^\times}

\newcommand{\lb}{\left(}
\newcommand{\rb}{\right)}
\newcommand{\lc}{\left\{}
\newcommand{\rc}{\right\}}
\newcommand{\relmid}{\mathrel{}\middle|\mathrel{}}
\DeclarePairedDelimiter{\abs}{\lvert}{\rvert}

\theoremstyle{plain}
\newtheorem{theorem}{Theorem}
\theoremstyle{definition}
\newtheorem{condition}{Condition}

\title[Unimodal singularities]{A new characterization of exceptional unimodal singularities}
\author[Y.~Katase]{Yuriko Katase}
\address{
Graduate School of Mathematical Sciences,
The University of Tokyo,
3-8-1 Komaba,
Meguro-ku,
Tokyo,
153-8914,
Japan.}
\email{ykatase@ms.u-tokyo.ac.jp}

\begin{document}

\maketitle

\begin{abstract}
Motivated by mirror symmetry for weighted projective spaces,
we give a new characterization of
exceptional unimodal singularities.
\end{abstract}

\section{introduction}

Exceptional unimodal singularities are introduced by Arnold
\cite{MR0467795}
as hypersurface singularities of modality one
which do not come in infinite family.
There are exactly 14 of them
which are defined by weighted homogeneous polynomials
as shown in Table \ref{tb:EUS}.
Dolgachev \cite{MR0345974} found another characterization
of exceptional unimodal singularities
as triangle singularities
which are hypersurfaces.
The triple $(\delta_1, \delta_2, \delta_3)$ of integers
specifying the triangle singularity
is called the \emph{Dolgachev number}.

The Milnor lattices of exceptional unimodal singularities
are computed by Gabrielov \cite{MR0367274}.
The Coxeter--Dynkin diagram
of the Milnor lattice
with respect to a distinguished basis of vanishing cycles
is specified by a triple $(\gamma_1, \gamma_2, \gamma_3)$
of integers
called the \emph{Gabrielov number}.
In Table \ref{tb:EUS},
one can see
that exceptional unimodal singularities
come in pairs in such a way that
the Dolgachev number and the Gabrielov number are interchanged.
This fact is discovered by Arnold
and given the name \emph{strange duality}.
Pinkham \cite{MR0429876}
and Dolgachev and Nikulin
\cite{MR728992,MR525944}
gave an interpretation of strange duality
as the exchange of the transcendental lattices
and algebraic lattices of K3 surfaces.
This interpretation can be considered as a precursor
of mirror symmetry for K3 surfaces
\cite{MR1416354,MR1420220}.
See e.g.~\cite{MR3044454,1407.1566,1806.04345}
and references therein
for the relation between exceptional unimodal singularities
and mirror symmetry for K3 surfaces.

\begin{table}[t]
\[
\begin{array}{cccccc}
  \toprule
 \text{name} & \text{normal form} & \text{weight system} &
 \text{D\#} & \text{G\#} & \text{dual} \\
  \midrule
 E_{12} & x^3 + y^7 + z^2
  & (6,14,21;42) & (2,3,7) & (2,3,7) & E_{12} \\
 E_{13} & x^3 + x y^5 + z^2
  & (4,10,15;30) & (2,4,5) & (2,3,8) & Z_{11} \\
 E_{14} & x^3 + y^8 + z^2
  & (3,8,12;24) & (3,3,4) & (2,3,9) & Q_{10} \\
 Z_{11} & x^3 y + y^5 + z^2
  & (6,8,15;30) & (2,3,8) & (2,4,5) & E_{13} \\
 Z_{12} & x^3 y + x y^4 + z^2
  & (4,6,11;22) & (2,4,6) & (2,4,6) & Z_{12} \\
 Z_{13} & x^3 y + y^6 + z^2
  & (3,5,9,18) & (3,3,5) & (2,4,7) & Q_{11} \\
 W_{12} & x^4 + y^5 + z^2
  & (4,5,10;20) & (2,5,5) & (2,5,5) & W_{12} \\
 W_{13} & x^4 + x y^4 + z^2
  & (3,4,8,16) & (3,4,4) & (2,5,6) & S_{11} \\
 Q_{10} & x^3 + y^4 + y z^2
  & (6,8,9;24) & (2,3,9) & (3,3,4) & E_{14} \\
 Q_{11} & x^3 + y^2 z + x z^3
  & (4,6,7;18) & (2,4,7) & (3,3,5) & Z_{13} \\
 Q_{12} & x^3 + y^5 + y z^2
  & (3,5,6;15) & (3,3,6) & (3,3,6) & Q_{12} \\
 S_{11} & x^4 + y^2 z + x z^2
  & (4,5,6,16) & (2,5,6) & (3,4,4) & W_{13} \\
 S_{12} & x^3 y + y^2 z + x z^2
  & (3,4,5;13) & (3,4,5) & (3,4,5) & S_{12} \\
 U_{12} & x^3 + y^3 + z^4
  & (3,4,4;12) & (4,4,4) & (4,4,4) & U_{12} \\
  \bottomrule
\end{array}
\]
\caption{14 exceptional unimodal singularities}
\label{tb:EUS}
\end{table}

A \emph{weight system} is a quadruple $(a_1,a_2,a_3;h)$ of positive integers.
The triple $a = (a_1,a_2,a_3)$ is called the \emph{weight},
and the integer $h$ is called the \emph{degree}.
A polynomial
$
 \sum_{i=(i_1,i_2,i_3) \in \bN^3} \alpha_i x^{i_1} y^{i_2} z^{i_3}
  \in \bC[x,y,z]
$
is \emph{weighted homogeneous of weight $a$ and degree $h$}
if $\sum_{k=1}^3 i_k a_k = h$ for any $i \in \bN^3$ with $\alpha_i \ne 0$.
We set
$
 a_0 = h - a_1 - a_2 - a_3.
$
It is known
by Reid (unpublished)
and Yonemura \cite{MR1066667}
that there are exactly 95 weights $a$
such that 
the minimal model
of a general anticanonical hypersurface
in the weighted projective space $\bP(a_0, a_1, a_2, a_3)$
is a K3 surface.
The weight systems in Table \ref{tb:EUS}
associated with
exceptional unimodal singularities
are on the Reid--Yonemura list.
Since all of them satisfy the condition $a_0 = 1$,
we will always assume this condition in this paper,
and write
$
 W_a = (a_1,a_2,a_3;a_1+a_2+a_3+1).
$
There are 41 weights with $a_0 = 1$
on the Reid--Yonemura list.

The mirror of
$
 \bP(a_0, a_1, a_2, a_3)
$
is given by the regular function
$
 F_a = u_0 + u_1 + u_2 + u_3
$
on
$
 M_q =
  \{ (u_0,u_1,u_2,u_3) \in \lb \bCx \rb^4
   \mid
    u_0^{a_0} u_1^{a_1} u_2^{a_2} u_3^{a_3} = q \},
$
whose fiber $F_a^{-1}(t)$ is the mirror of an anticanonical hypersurface
in $\bP(a_0, a_1, a_2, a_3)$
\cite{MR1403947}.
We set the K\"ahler parameter $q$ to 1 for simplicity.
Since $a_0 = 1$,
the fiber $F_a^{-1}(t)$ can be described as
\begin{align}
 F_a^{-1}(t)
  &\cong \lc (u_1,u_2,u_3) \in \lb \bCx \rb^3 \relmid
 u_1 + u_2 + u_3 + \frac{1}{u_1^{a_1} u_2^{a_2} u_3^{a_3}} = t \rc \\
  &= \lc (u_1,u_2,u_3) \in \lb \bCx \rb^3 \relmid
 u_1^{a_1+1} u_2^{a_2} u_3^{a_3}
  + u_1^{a_1} u_2^{a_2+1} u_3^{a_3}
  + u_1^{a_1} u_2^{a_2} u_3^{a_3+1} + 1 
  = t u_1^{a_1} u_2^{a_2} u_3^{a_3} \rc.
\end{align}
We define polynomial maps $G_a, H_a \colon \bC^3 \to \bC$ by
\begin{align}
 G_a(u_1,u_2,u_3) &= u_1^{a_1+1} u_2^{a_2} u_3^{a_3}
  + u_1^{a_1} u_2^{a_2+1} u_3^{a_3}
  + u_1^{a_1} u_2^{a_2} u_3^{a_3+1}, \\
 H_a(u_1,u_2,u_3) &= u_1^{a_1} u_2^{a_2} u_3^{a_3},
\end{align}
so that $F_a^{-1}(t)$ is defined by
\begin{align}
 G_a(u_1,u_2,u_3) + 1 = t H_a(u_1,u_2,u_3).
\end{align}
An integer matrix $C = \lb c_{ij} \rb_{i,j=1}^3$
with non-negative entries
defines a polynomial map
$
 f_C \colon \bC^3 \to \bC
$
by
\begin{align}
  f_C(x,y,z) = \sum_{i=1}^3 x^{c_{i1}} y^{c_{i2}} z^{c_{i3}}.
\end{align}
Similarly,
an integer matrix $D = \lb d_{ij} \rb_{i,j=1}^3$
with non-negative entries
defines a polynomial map
$
 \varphi_D \colon \bC^3 \to \bC^3
$
by
\begin{align}
 \varphi_D(u_1,u_2,u_3)
  = \lb u_1^{d_{i1}} u_2^{d_{i2}} u_3^{d_{i3}} \rb_{i=1}^3,
\end{align}
which restricts
to an isomorphism of tori $(\bCx)^3$
if $\abs{\det D} = 1$.
We also define a polynomial map $h \colon \bC^3 \to \bC$ by
\begin{align}
 h(x,y,z) = x y z.
\end{align}
Consider the following condition on the weight $a$:

\begin{condition} \label{cd:dual}
There exist $3 \times 3$ integer matrices
$C = (c_{ij})_{i,j=1}^3$
and
$D = (d_{ij})_{i,j=1}^3$
with non-negative entries
satisfying
\begin{enumerate}[(i)]
 \item \label{it:D}
$
 \abs{\det D} = 1,
$
 \item \label{it:G}
$
 G_a = f_C \circ \varphi_D,
$
 \item \label{it:H}
$
 H_a = h \circ \varphi_D,
$
and
 \item \label{it:f}
$f_C$
has an isolated critical point at the origin.
\end{enumerate}
\end{condition}

If Condition \ref{cd:dual} is satisfied,
then $F_a^{-1}(t)$ can be compactified to a quasi-smooth hypersurface
\begin{align}
 f_C(x,y,z) + w^{\hv} = t x y z w
\end{align}
of Dwork type in $\bP(1,\av_1,\av_2,\av_3)$,
where
\begin{align} \label{eq:av}
 \av_i = \sum_{j=1}^3 d_{ij}
\end{align}
and
$
 \hv = 1 + \sum_{i=1}^3 \av_i.
$
Note that Conditions (\ref{cd:dual}.\ref{it:G})
and (\ref{cd:dual}.\ref{it:H}) can be written as
\begin{align} \label{eq:BCD}
 B_a = C D
\end{align}
and
\begin{align} \label{eq:decomp_a}
\begin{pmatrix}
 a_1 & a_2 & a_3
\end{pmatrix}
 =
\begin{pmatrix}
 1 & 1 & 1
\end{pmatrix}
 D
\end{align}
respectively,
where
\begin{align}
 B_a =
\begin{pmatrix}
 a_1 + 1 & a_2 & a_3 \\
 a_1 & a_2 + 1 & a_3 \\
 a_1 & a_2 & a_3 + 1
\end{pmatrix}.
\end{align}
It follows from \eqref{eq:BCD} and
Condition (\ref{cd:dual}.\ref{it:D}) that
\begin{align}
 \abs{\det C} = \det B_a = a_1+a_2+a_3+1.
\end{align}

The main result in this paper
is the following characterization
of exceptional unimodal singularities:

\begin{theorem} \label{th:main}
A weight $(a_0,a_1,a_2,a_3)$ with $a_0=1$
on the Reid--Yonemura list comes from an exceptional unimodal singularity
if and only if $a = (a_1,a_2,a_3)$ satisfies Condition \ref{cd:dual}.
If this is the case, then $C$ is uniquely determined by $a$
up to permutation of rows and columns,
and $f_C$ is a defining polynomial
of the strange dual singularity.
\end{theorem}

For each $a$ among 41 on the Yonemura--Reid list,
there are finitely many integer matrices $D$
with non-negative entries
satisfying \eqref{eq:decomp_a}.
Each such $D$ determines $C$ by \eqref{eq:BCD},
and the complete list of $C$ and $D$
such that $\abs{\det D} = 1$,
up to permutation of rows and columns,
is shown in Table \ref{tb:semi-dual}.
If this is the case,
then we say that the weight $\av$
defined by \eqref{eq:av} is \emph{semi-dual}
to the weight $a$.
It follows from Table \ref{tb:semi-dual}
that semi-duality is reflexive;
$\av$ is semi-dual to $a$
if and only if $a$ is semi-dual to $\av$.
The proof Theorem \ref{th:main} is given by
testing if each $C$ satisfies Condition \ref{cd:dual}.(\ref{it:f}).

\begin{longtable}
{ccccc}
\toprule
No. & $a$ & $C$ & $D$ & semi-dual \\
\midrule
1 & $(1,1,1)$ 
&$\begin{pmatrix}
             $2$ & $1$ & $1$ \\
             $1$ & $2$ & $1$ \\
             $1$ & $1$ & $2$
            \end{pmatrix}$
&$\begin{pmatrix}
             $1$ & $0$ & $0$ \\
             $0$ & $1$ & $0$ \\
             $0$ & $0$ & $1$
            \end{pmatrix}$
& 1 \\
3 & $(2,2,1)$ 
& $\begin{pmatrix}
             $1$ & $1$ & $2$ \\
             $2$ & $1$ & $0$ \\
             $1$ & $2$ & $0$
            \end{pmatrix}$
&
            $\begin{pmatrix}
             $1$ & $1$ & $0$ \\
             $1$ & $0$ & $1$ \\
             $1$ & $0$ & $0$
            \end{pmatrix}$
&5\\
{} & ${}$ 
& $\begin{pmatrix}
             $2$ & $1$ & $0$ \\
             $1$ & $1$ & $2$ \\
             $0$ & $2$ & $2$
            \end{pmatrix}$
&
            $\begin{pmatrix}
             $1$ & $1$ & $0$ \\
             $1$ & $0$ & $1$ \\
             $0$ & $1$ & $0$
            \end{pmatrix}$
&3\\
4 & $(4,4,3)$ \
& $\begin{pmatrix}
             $2$ & $1$ & $0$ \\
             $1$ & $2$ & $0$ \\
             $0$ & $0$ & $4$
            \end{pmatrix}$
&
            $\begin{pmatrix}
             $2$ & $1$ & $1$ \\
             $1$ & $2$ & $1$ \\
             $1$ & $1$ & $1$
            \end{pmatrix}$
&4 (dual)\\
5 & $(3,1,1)$ 
& $\begin{pmatrix}
             $1$ & $2$ & $1$ \\
             $1$ & $1$ & $2$ \\
             $2$ & $0$ & $0$
            \end{pmatrix}$
&
            $\begin{pmatrix}
             $1$ & $1$ & $1$ \\
             $1$ & $0$ & $0$ \\
             $0$ & $1$ & $0$
            \end{pmatrix}$
&3\\
6 & $(5,2,2)$ 
& $\begin{pmatrix}
             $2$ & $0$ & $0$ \\
             $0$ & $3$ & $2$ \\
             $0$ & $2$ & $3$
            \end{pmatrix}$
&
            $\begin{pmatrix}
             $3$ & $1$ & $1$ \\
             $1$ & $1$ & $0$ \\
             $1$ & $0$ & $1$
            \end{pmatrix}$
&6\\
7 & $(4,2,1)$
& $\begin{pmatrix}
             $1$ & $1$ & $2$ \\
             $1$ & $2$ & $0$ \\
             $2$ & $0$ & $0$
            \end{pmatrix}$
 &
            $\begin{pmatrix}
             $2$ & $1$ & $1$ \\
             $1$ & $1$ & $0$ \\
             $1$ & $0$ & $0$
            \end{pmatrix}$
&7\\
{} & ${}$
& $\begin{pmatrix}
             $2$ & $0$ & $0$ \\
             $0$ & $3$ & $2$ \\
             $1$ & $1$ & $2$
            \end{pmatrix}$
 &
            $\begin{pmatrix}
             $2$ & $1$ & $1$ \\
             $1$ & $1$ & $0$ \\
             $0$ & $0$ & $1$
            \end{pmatrix}$
&19\\
8 & $(6,3,2)$ 
& $\begin{pmatrix}
             $1$ & $2$ & $0$ \\
             $1$ & $0$ & $3$ \\
             $2$ & $0$ & $0$
            \end{pmatrix}$
&
            $\begin{pmatrix}
             $3$ & $2$ & $1$ \\
             $2$ & $1$ & $0$ \\
             $1$ & $1$ & $0$
            \end{pmatrix}$
&10\\
{} & ${}$ 
& $\begin{pmatrix}
             $2$ & $0$ & $0$ \\
             $0$ & $2$ & $3$ \\
             $1$ & $2$ & $0$
            \end{pmatrix}$
&
            $\begin{pmatrix}
             $3$ & $2$ & $1$ \\
             $1$ & $1$ & $1$ \\
             $1$ & $1$ & $0$
            \end{pmatrix}$
&24\\
9 & $(10,5,4)$ 
& $\begin{pmatrix}
             $1$ & $2$ & $0$ \\
             $2$ & $0$ & $0$ \\
             $0$ & $0$ & $5$
            \end{pmatrix}$
&
            $\begin{pmatrix}
             $5$ & $3$ & $2$ \\
             $3$ & $1$ & $1$ \\
             $2$ & $1$ & $1$
            \end{pmatrix}$
&9 (dual)\\
10 & $(6,4,1)$ 
& $\begin{pmatrix}
             $1$ & $1$ & $2$ \\
             $2$ & $0$ & $0$ \\
             $0$ & $3$ & $0$
            \end{pmatrix}$
&
            $\begin{pmatrix}
             $3$ & $2$ & $1$ \\
             $2$ & $1$ & $1$ \\
             $1$ & $0$ & $0$
            \end{pmatrix}$
&8\\
12 & $(9,6,2)$ 
& $\begin{pmatrix}
             $2$ & $0$ & $0$ \\
             $0$ & $2$ & $3$ \\
             $0$ & $3$ & $0$
            \end{pmatrix}$
&
            $\begin{pmatrix}
             $5$ & $3$ & $1$ \\
             $3$ & $2$ & $1$ \\
             $1$ & $1$ & $0$
            \end{pmatrix}$
&12\\
13 & $(12,8,3)$ 
&$\begin{pmatrix}
             $2$ & $0$ & $0$ \\
             $0$ & $3$ & $0$ \\
             $1$ & $0$ & $4$
            \end{pmatrix}$
&  
            $\begin{pmatrix}
             $5$ & $4$ & $3$ \\
             $3$ & $3$ & $2$ \\
             $1$ & $1$ & $1$
            \end{pmatrix}$
& 20 (dual) \\
14 & $(21,14,6)$ 
&$\begin{pmatrix}
             $2$ & $0$ & $0$ \\
             $0$ & $3$ & $0$ \\
             $0$ & $0$ & $7$
            \end{pmatrix}$
& 
            $\begin{pmatrix}
             $11$ & $7$ & $3$ \\
             $7$ & $5$ & $2$ \\
             $3$ & $2$ & $1$
            \end{pmatrix}$
& 14 (dual)\\
18 & $(3,3,2)$ 
& $\begin{pmatrix}
             $1$ & $2$ & $0$ \\
             $1$ & $0$ & $3$ \\
             $2$ & $1$ & $0$
            \end{pmatrix}$
&
            $\begin{pmatrix}
             $1$ & $1$ & $1$ \\
             $2$ & $1$ & $0$ \\
             $1$ & $1$ & $0$
            \end{pmatrix}$
&25\\
19 & $(3,2,2)$ 
& $\begin{pmatrix}
             $2$ & $0$ & $1$ \\
             $0$ & $3$ & $1$ \\
             $0$ & $2$ & $2$
            \end{pmatrix}$
&
            $\begin{pmatrix}
             $2$ & $1$ & $0$ \\
             $1$ & $1$ & $0$ \\
             $1$ & $0$ & $1$
            \end{pmatrix}$
&7\\
{} & ${}$
& $\begin{pmatrix}
             $0$ & $2$ & $2$ \\
             $2$ & $0$ & $1$ \\
             $2$ & $1$ & $0$
            \end{pmatrix}$
 &
            $\begin{pmatrix}
             $1$ & $1$ & $1$ \\
             $1$ & $0$ & $1$ \\
             $1$ & $1$ & $0$
            \end{pmatrix}$
&19\\
20 & $(9,8,6)$
&$\begin{pmatrix}
             $2$ & $0$ & $1$ \\
             $0$ & $3$ & $0$ \\
             $0$ & $0$ & $4$
            \end{pmatrix}$
 &
            $\begin{pmatrix}
             $5$ & $3$ & $1$ \\
             $4$ & $3$ & $1$ \\
             $3$ & $2$ & $1$
            \end{pmatrix}$ 
& 13 (dual)\\
21 & $(2,1,1)$ 
& $\begin{pmatrix}
             $1$ & $2$ & $1$ \\
             $2$ & $0$ & $1$ \\
             $1$ & $1$ & $2$
            \end{pmatrix}$
&
            $\begin{pmatrix}
             $1$ & $1$ & $0$ \\
             $1$ & $0$ & $0$ \\
             $0$ & $0$ & $1$
            \end{pmatrix}$
&21\\
22 & $(6,5,3)$ 
&$\begin{pmatrix}
             $2$ & $0$ & $1$ \\
             $0$ & $3$ & $0$ \\
             $1$ & $0$ & $3$
            \end{pmatrix}$
& 
            $\begin{pmatrix}
             $3$ & $2$ & $1$ \\
             $2$ & $2$ & $1$ \\
             $1$ & $1$ & $1$
            \end{pmatrix}$
& 22 (dual)\\
24 & $(5,4,2)$ 
& $\begin{pmatrix}
             $2$ & $0$ & $1$ \\
             $0$ & $2$ & $2$ \\
             $0$ & $3$ & $0$
            \end{pmatrix}$
&
            $\begin{pmatrix}
             $3$ & $1$ & $1$ \\
             $2$ & $1$ & $1$ \\
             $1$ & $1$ & $0$
            \end{pmatrix}$
&8\\
25 & $(4,3,1)$ 
& $\begin{pmatrix}
             $2$ & $0$ & $1$ \\
             $1$ & $1$ & $2$ \\
             $0$ & $3$ & $0$
            \end{pmatrix}$
&
            $\begin{pmatrix}
             $2$ & $1$ & $1$ \\
             $1$ & $1$ & $1$ \\
             $0$ & $1$ & $0$
            \end{pmatrix}$
&18\\
28 & $(10,7,3)$ & -- & -- & -- \\
37 & $(8,4,3)$ 
&$\begin{pmatrix}
             $1$ & $2$ & $0$ \\
             $2$ & $0$ & $0$ \\
             $0$ & $1$ & $4$
            \end{pmatrix}$
&
            $\begin{pmatrix}
             $3$ & $3$ & $2$ \\
             $2$ & $1$ & $1$ \\
             $1$ & $1$ & $1$
            \end{pmatrix}$ 
& 58 (dual)\\
38 & $(15,8,6)$ 
&$\begin{pmatrix}
             $2$ & $0$ & $0$ \\
             $0$ & $3$ & $1$ \\
             $0$ & $0$ & $5$
            \end{pmatrix}$
& 
            $\begin{pmatrix}
             $8$ & $5$ & $2$ \\
             $4$ & $3$ & $1$ \\
             $3$ & $2$ & $1$
            \end{pmatrix}$ 
& 50 (dual)\\
39 & $(9,5,3)$ 
&$\begin{pmatrix}
             $2$ & $0$ & $0$ \\
             $0$ & $3$ & $1$ \\
             $1$ & $0$ & $3$
            \end{pmatrix}$
& 
            $\begin{pmatrix}
             $4$ & $3$ & $2$ \\
             $2$ & $2$ & $1$ \\
             $1$ & $1$ & $1$
            \end{pmatrix}$ 
&60 (dual)\\
40 & $(7,4,2)$ 
& $\begin{pmatrix}
             $2$ & $0$ & $0$ \\
             $0$ & $2$ & $3$ \\
             $0$ & $3$ & $1$
            \end{pmatrix}$
&
            $\begin{pmatrix}
             $4$ & $2$ & $1$ \\
             $2$ & $1$ & $1$ \\
             $1$ & $1$ & $0$
            \end{pmatrix}$
&40\\
42 & $(5,3,1)$ 
& $\begin{pmatrix}
             $1$ & $1$ & $2$ \\
             $2$ & $0$ & $0$ \\
             $0$ & $3$ & $1$
            \end{pmatrix}$
&
            $\begin{pmatrix}
             $2$ & $2$ & $1$ \\
             $1$ & $1$ & $1$ \\
             $1$ & $0$ & $0$
            \end{pmatrix}$
&63\\
44 & $(8,5,2)$ & -- & -- & -- \\
45 & $(14,9,4)$ & -- & -- & -- \\
50 & $(15,10,4)$ 
&$\begin{pmatrix}
             $2$ & $0$ & $0$ \\
             $0$ & $3$ & $0$ \\
             $0$ & $1$ & $5$
            \end{pmatrix}$
& 
            $\begin{pmatrix}
             $8$ & $4$ & $3$ \\
             $5$ & $3$ & $2$ \\
             $2$ & $1$ & $1$
            \end{pmatrix}$ 
&38 (dual)\\
51 & $(18,12,5)$ & -- & -- & -- \\
58 & $(6,5,4)$ 
&$\begin{pmatrix}
             $1$ & $2$ & $0$ \\
             $2$ & $0$ & $1$ \\
             $0$ & $0$ & $4$
            \end{pmatrix}$
&
            $\begin{pmatrix}
             $3$ & $2$ & $1$ \\
             $3$ & $1$ & $1$ \\
             $2$ & $1$ & $1$
            \end{pmatrix}$
&37 (dual)\\
59 & $(8,7,5)$ & -- & -- & -- \\
60 & $(7,6,4)$ 
&$\begin{pmatrix}
             $2$ & $0$ & $1$ \\
             $0$ & $3$ & $0$ \\
             $0$ & $1$ & $3$
            \end{pmatrix}$
&
            $\begin{pmatrix}
             $4$ & $2$ & $1$ \\
             $3$ & $2$ & $1$ \\
             $2$ & $1$ & $1$
            \end{pmatrix}$ 
&39 (dual)\\
63 & $(4,3,2)$
& $\begin{pmatrix}
             $1$ & $2$ & $0$ \\
             $1$ & $0$ & $3$ \\
             $2$ & $0$ & $1$
            \end{pmatrix}$
 &
            $\begin{pmatrix}
             $2$ & $1$ & $1$ \\
             $2$ & $1$ & $0$ \\
             $1$ & $1$ & $0$
            \end{pmatrix}$
&42\\
66 & $(3,2,1)$ 
& $\begin{pmatrix}
             $2$ & $0$ & $1$ \\
             $0$ & $3$ & $1$ \\
             $1$ & $1$ & $2$
            \end{pmatrix}$
&
            $\begin{pmatrix}
             $2$ & $1$ & $0$ \\
             $1$ & $1$ & $0$ \\
             $0$ & $0$ & $1$
            \end{pmatrix}$
&66\\
{} & ${}$
& $\begin{pmatrix}
             $1$ & $1$ & $2$ \\
             $1$ & $2$ & $0$ \\
             $2$ & $0$ & $1$
            \end{pmatrix}$
 &
            $\begin{pmatrix}
             $1$ & $1$ & $1$ \\
             $1$ & $1$ & $0$ \\
             $1$ & $0$ & $0$
            \end{pmatrix}$
&66\\
71 & $(7,4,3)$ & -- & -- & -- \\
72 & $(7,5,2)$ & -- & -- & -- \\
77 & $(13,7,5)$ & -- & -- & -- \\
78 & $(11,6,4)$ 
&$\begin{pmatrix}
             $2$ & $0$ & $0$ \\
             $0$ & $3$ & $1$ \\
             $0$ & $1$ & $4$
            \end{pmatrix}$
& 
           $\begin{pmatrix}
             $6$ & $3$ & $2$ \\
             $3$ & $2$ & $1$ \\
             $2$ & $1$ & $1$
            \end{pmatrix}$ 
&78 (dual)\\
82 & $(11,7,3)$ & -- & -- & -- \\
87 & $(5,4,3)$ 
&$\begin{pmatrix}
             $1$ & $2$ & $0$ \\
             $2$ & $0$ & $1$ \\
             $0$ & $1$ & $3$
            \end{pmatrix}$
& 
            $\begin{pmatrix}
             $2$ & $2$ & $1$ \\
             $2$ & $1$ & $1$ \\
             $1$ & $1$ & $1$
            \end{pmatrix}$ 
&87 (dual)\\
89 & $(5,3,2)$ & -- & -- & -- \\
\bottomrule\\
\caption{Semi-duality of weights}
\label{tb:semi-dual}
\end{longtable}

Note that a normal form of a weighted homogeneous
exceptional unimodal singularity is not unique in general,
and Theorem \ref{th:main} allows us to fix one uniquely.
For example,
the defining equation for
the $W_{12}$-singularity can be written either as
$x^5 + y^4 + z^2$ or $x^5 + y^2 z + z^2$,
and only the latter comes from Theorem \ref{th:main}.

Recall from \cite{MR2426805} that weight systems
$W_a = (a_1,a_2,a_3;h)$
and
$W_{\av} = \lb \av_1,\av_2,\av_3;\hv \rb$
are said to be \emph{Kobayashi dual}
if there is a $3\times 3$ integer matrix $C$
with non-negative entries
satisfying the \emph{weighted magic square condition}
\begin{align}
 C
\begin{pmatrix}
 a_1 \\ a_2 \\ a_3
\end{pmatrix}
 =
\begin{pmatrix}
 h \\ h \\ h
\end{pmatrix}
 \quad \text{and} \quad
\begin{pmatrix}
 \av_1 & \av_2 & \av_3
\end{pmatrix}
 C =
\begin{pmatrix}
 \hv & \hv & \hv
\end{pmatrix}
\end{align}
and the \emph{primitivity}
\begin{align}
 \abs{\det C}=h=\hv.
\end{align}
Kobayashi duality is a generalization of strange duality
to weights which may not come from exceptional unimodal singularities
\cite{MR2426805,MR2278769}.
The matrices $C$ appearing in Table \ref{tb:semi-dual}
are primitive weighted magic squares,
so that our semi-duality is a special case of Kobayashi duality.

\ \\
\emph{Acknowledgements.}
The author thanks her advisor Kazushi Ueda for guidance and encouragement.

\bibliographystyle{amsalpha}
\bibliography{bibs}

\end{document}